\magnification=1200
\def\n{\hfill\break}
\centerline{\bf Grassmannian and Elliptic Operators}
\vskip 0.3in
\centerline{Leonid Friedlander\footnote*{Partially supported by NSF}}
\centerline{University of Arizona}
\centerline{friedlan@math.arizona.edu}
\centerline{Albert Schwarz$^*$}
\centerline{University of California, Davis}
\centerline{schwarz@math.ucdavis.edu}
\vskip 0.3in
\n
{\bf Abstract.} {\it A multi-dimensional analogue of the Krichever 
construction is discussed.}
\vskip 0.3in
\centerline{\bf 1. Introduction}

Infinite-dimensional Grassmannian (the Sato Grassmannian) plays an important
role in mathematical physics, as well as in other branches of mathematics. 
It was conjectured, in particular, that it constitutes a natural
framework for non-perturbative formulation of the string theory [1--5].
This conjecture is based first of all on the remark that moduli
spaces of algebraic curves of all genera are embedded in the Grassmannian by
means of so called Krichever construction. New evidence in favor of
this conjecture was presented recently in [6].

The development of the string theory led to a conclusion that strings
should not be considered as fundamental objects. They should appear
in broader theory on equal footing with  their multidimensional
analogues -- membranes. If we believe that this broader theory
can be formulated in terms of infinite-dimensional Grassmannians,
we can conjecture that membranes are related to the Grassmannian
by means of multidimensional generalization of the
Krichever construction. The main goal of this paper is to describe
such a generalization.

Using the well-known relation between the $\bar\partial$-operator
and the Dirac operator we can formulate the Krichever construction
in terms of the Dirac operator; its multidimensional generalization
also can be formulated in this way.

Namely, we shall consider a Riemannian manifold $M$ with boundary
$\Gamma$, a Hermitian vector bundle $E$ over $M$, and a connection
in this bundle. These data permit us to construct corresponding
Dirac operator $A$ (if some topological conditions are satisfied).
In physicist's terms, we consider the Dirac operator in an external
gauge field. This operator acts in the space of smooth sections
of the vector bundle $S\otimes E$ where $S$ stands for the spinor bundle
over $M$. One can introduce a natural Hermitian inner product
in this space of sections, and one can construct a Hilbert space as the
completion of the space of smooth sections with respect  to this inner product.
The Dirac operator becomes a self-adjoint operator in this Hilbert space.
One can also introduce the boundary Hilbert space ${\cal H}_\Gamma$,
the completion of the space of smooth sections of $S\otimes E$ over
$\Gamma$. We define the subspace $H_+^A$ of ${\cal H}_\Gamma$ as the closure
of the space of all sections of $S\otimes E$ over $\Gamma$ that can be extended
to smooth solutions of the Dirac equation $Af=0$ over $M$.
Using the decomposition of ${\cal H}_\Gamma$ in the direct sum of $H_+^A$
and its orthogonal complement, $H_-^A$, we construct the Segal--Wilson version
of the Sato Grassmannian $Gr$ in the standard way: a linear subspace
$V\subset {\cal H}_\Gamma$ belongs to $Gr$ if the projection
$\pi_+:V\to H_+^A$ is a Fredholm operator and the projection
$\pi_-:V\to H_-^A$ is a compact operator. Here $\pi_\pm$ is the
orthogonal projection onto $H_\pm^A$.

Let us replace now the manifold $M$ by a manifold $\tilde M$ that has
the same boundary $\Gamma$, and, instead of the Hermitian vector bundle
$p_E: E\to M$, let us take a Hermitian vector bundle 
$p_{\tilde E}: \tilde E\to\tilde M$.
The new objects are taken in such a way that they coincide with old
ones on $\Gamma$. More accurately,\n
({\it i}) Riemannian metrics on $\Gamma$ that are induced by 
Riemannian metrics on $M$ and $\tilde M$ coincide;\n
({\it ii}) there exist a neighborhood $U$ of $\Gamma$ in $M$,
a neighborhood $\tilde U$ of $\Gamma$ in $\tilde M$, and an isomorphism
$\Psi:\tilde E_{\tilde U}\to E_U$ (here $E_U$ is the restriction
of $E$ to $U$,...) such that\n
(a) $p_{\tilde E}=p_E\circ\Psi$ over $\Gamma$;\n
(b) $\Psi$ is an isometry over $\Gamma$ (when restricted to $\tilde E_\Gamma$);
\n
(c) $\Psi^{-1}\nabla\Psi -\tilde\nabla$ is a differential operator
of order $0$ on $\Gamma$. Here $\nabla$ and $\tilde\nabla$
are connections on $E$ and $\tilde E$, respectively.

One can define new Dirac operator $\tilde A$, and this operator gives rise
to the subspace $H_+^{\tilde A}$ of ${\cal H}_\Gamma$. We will prove that
$H_+^{\tilde A}$ belongs to the Grassmannian. More precisely,
$\pi_+: H_+^{\tilde A}\to H_+^A$ is a Fredholm operator of index $0$,
and the operator $\pi_-: H_+^{\tilde A}\to H_-^A$ belongs to the
Schatten ideal $\Sigma_p$ for $p>\dim M-1$.

The spinor bundle $S$ over $M$ can be decomposed into the direct sum
of the bundle of left spinors, $^LS$, and the bundle of right spinors,
$^RS$. The Dirac operator takes a section of $^LS\otimes E$
to a section of $^RS\otimes E$ and vice versa. 
We denote by $^{L(R)}{\cal H}_\Gamma$ the space of smooth sections of 
$^{L(R)}S\otimes E$, completed with respect to the inner product
that was discussed earlier. This decomposition
gives us decompositions of the subspaces $H_\pm^A$ into direct sums
${H_\pm^A=^L\!\!H_\pm^A\oplus ^R\!\!H_\pm^A}$. Clearly,
$^L{\cal H}_\Gamma=^L\!\!\!H_+^A\oplus ^L\!\!H_-^A$.  In the same way as 
we defined
the Grassmannian of subspaces in ${\cal H}_\Gamma$, one defines the 
Grassmannian
in $^L{\cal H}_\Gamma$, and $^LH_+^{\tilde A}$ belongs to this Grassmannian.
However, the index of the Fredholm operator 
$\pi_+:^L\!\!H_+^{\tilde A}\to ^L\!\!H_+^A$ is not necessarily zero.

In the case when $\dim M=2$, one can identify the spaces $^LH_+^{\tilde A}$
with the points in $Gr$ that can be obtained from the standard Krichever
construction. Notice that in this case the operators
$\pi_-:^L\!\!H_+^{\tilde A}\to ^L\!\!H_-^A$ belong, in particular, to the
Schatten ideal $\Sigma_2$, i.e. they are Hilbert--Schmidt operators. 
The corresponding points in the Grassmannian have an interpretation in
terms of the fermionic Fock space. As we mentioned earlier, in the case
when the dimension of $M$ is arbitrary, these operators belong
to the Schatten ideal $\Sigma_p$, with $p>\dim M -1$.
The fermionic interpretation of the corresponding points in the Grassmannian
was analyzed in [7].

The points in the Grassmannian that can be obtained  by means of the 
Krichever
construction have ``large stabilizers'' in appropriate groups acting
on the Grassmannian, and they can be characterized by this property.
This fact plays an important role in [6]. It would be interesting to obtain
similar results for the multidimensional analogue of the Krichever 
construction. We have only tentative results in this direction.

In the above statements, one can replace the multi-dimensional
Dirac operator by an arbitrary elliptic differential
operator. We will prove our main result in this generality.
The proof uses standard technique of the theory of
pseudodifferential operators. First, the proof will be carried out
under an additional assumption that the Agmon--Seeley condition
is satisfied. Then we will show that this assumption can be removed.
In the case of the Dirac operator,
one can use the considerations of [8] for deriving our results.
Our general proof is based on the technique of [9], and it follows [9]
rather closely.

The authors wish to thank G.Henkin, M.Kontsevich, and M.Shubin
for discussions.

It is a privilege for us to publish our paper in this volume.
J.Stasheff made outstanding contribution to mathematics.
The algebraic structures that he introduced and studied
play an important role in modern mathematical physics.
We dedicate our paper to J.Stasheff.

\centerline{\bf 2. A manifold without boundary}
Let $E\to M$ be a vector bundle of rank $r$ over a compact, orientable, 
closed manifold $M$ of dimension $n$, and let $\Gamma$ be a hypersurface
in $M$ that divides $M$ into two connected components, $M_+$ and $M_-$.
We denote by $E_\Gamma$ the restriction of $E$ to $\Gamma$, and
$E_\pm$ are restrictions of $E$ to $M_\pm$. One can find a function,
$x_n$, defined in a neighborhood of $\Gamma$ in $M$ such that
$x_n=0$ on $\Gamma$, $\pm x_n>0$ in $M_\pm$, and $dx_n\not= 0$.
When local coordinates in a neighborhood of a point from $\Gamma$ are used,
they will be always assumed to be of the form $(x',x_n)$ where
$x'=(x_1,\ldots x_{n-1})$ are local coordinates  on $\Gamma$.
We introduce a connection $\nabla$ on the restriction of $E$ to
a neighborhood of $\Gamma$, and, with a slight abuse of notations,
we will denote by $\partial_n$, or by $\partial/\partial x_n$,
the covariant derivative with respect to $x_n$.

Let $A$ be an elliptic pseudo-differential operator of order $k>0$
that acts on sections of $E$. Throughout this section,
we will assume that\n
({\it i}) the operator $A$ satisfies the Agmon--Seeley condition:
there exists an angle $\{z: |\arg z-\theta |<\epsilon\}$
in the complex plane that is free from eigenvalues of the principal
symbol of $A$;\n
({\it ii}) in a neighborhood of $\Gamma$, the operator $A$ is differential;\n
({\it iii}) the restriction of $Au$ to $M_\pm$ depends on the restriction
of $u$ to $M_\pm$ only.\n
For technical reasons, we also assume
that $0$ is a regular point of $A$. This assumption is not essential.
We will indicate, what changes should be made in the case when the operator 
$A$ is not invertible.
A neighborhood of $\Gamma$ in $M$ is diffeomorphic to
$\Gamma \times (-1,1)$. We fix this diffeomorphism once and forever;
$x_n$ is the coordinate along the interval $(-1,1)$.
Assume that in this neighborhood our operator $A$ is differential, and
it can be written as
$$A=A_k\partial_n^k+\cdots +A_1\partial_n+A_0$$
where $A_q$ is a differential operator of order $k-q$
that contains the tangential derivatives only. In particular
$A_k$ is a differential operator of order $0$, that is a smooth
family of endomorphisms of fibers of $E$. Ellipticity of
the operator $A$ implies that the endomorphisms
$A_k$ are non-degenerate. In fact, if one denotes by $\xi'$
the set of dual variables to $x'$, and by $\xi_n$ the variable dual to
$x_n$ then the principal symbol of $A$, when evaluated at a point
$(x;\xi'=0,\xi_n)$, equals $A_k(x)(i\xi_n)^k$. It should be non-degenerate
when $\xi_n\not= 0$.

We introduce two subspaces,
${\cal L}_\pm$, of the space
$${\cal L}=\underbrace{C^\infty(\Gamma,E_\Gamma)\oplus\cdots\oplus
  C^\infty(\Gamma,E_\Gamma)}_{k\rm\;times}$$
of sections of the vector bundle
$${\cal E}=\underbrace{E_\Gamma\oplus\cdots\oplus
E_\Gamma}_{k\rm\;times}$$
in the following way
$${\cal L}_\pm = \{(\phi_0,\ldots,\phi_{k-1}): \phi_j=\partial_n^j u \quad
\hbox{\rm where}\quad Au=0 \quad\hbox{\rm in}\quad M_\pm\}.$$

{\bf Proposition 1.} {\sl Under all assumptions made above,}
$${\cal L}_+\cap{\cal L}_-=\{0\};\quad {\cal L}_+ +{\cal L}_-={\cal L}.$$

{\bf Proof.} The first statement is almost obvious. Let 
$\phi=(\phi_0,\ldots,\phi_{k-1})\in {\cal L}_+\cap{\cal L}_-$. Then there exist
sections $u_\pm$ of $E_\pm$ such that $Au_\pm =0$ and
$\partial_n^j u_\pm =\phi_j$,
$j=0,\ldots , k-1$, on $\Gamma$. Then, from the equation $Au_\pm=0$,
it follows that all partial derivatives of $u_+$ and $u_-$ agree on
$\Gamma$ (we have used ellipticity of the operator $A$). Therefore
the section $u$ of $E$ defined as $u_+$ over $M_+$ and as $u_-$ over
$M_-$ is smooth; so it satisfies $Au=0$. Because $0$ does not
belong to the spectrum of $A$, we conclude that $u=0$, and $\phi=0$.

We proceed now to proving the second statement of the Proposition.
Let $u_\pm$ be a section of $E_\pm\to M_\pm$, smooth up to $\Gamma$.
By $u_\pm^0$ we denote the section of $E\to M$ that equals $u$ over
$M_\pm$, and that equals $0$ over $M_\mp$. Then, for every positive 
integer number $q$,
$$\partial_n^q(u_\pm^0)=(\partial_n^q u_\pm)^0 \pm
  \sum_{p=0}^{q-1}\phi_{q-p-1}^\pm \delta^{(p)}(x_n)$$
where $\phi_j^\pm$ is the restriction of $\partial_n^j u_\pm$ to $\Gamma$.
We will denote by ${\cal A}_q$ restrictions of operators $A_q$
to $\Gamma$. They are differential operators acting
on sections of the vector bundle $E_\Gamma\to\Gamma$. One has
$$\eqalignno{A(u_\pm^0) & = (Au_\pm)^0 \pm \sum_{q=1}^k{\cal A}_q
 \sum_{p=0}^{q-1}\phi_{q-p-1}^\pm \delta^{(p)}(x_n)\cr
  & = (Au_\pm)^0 \pm \sum_{p=0}^{k-1}\biggl( \sum_{q=p+1}^k
 {\cal A}_q\phi_{q-p-1}^\pm \biggr) \delta^{(p)}(x_n)\cr
  & = (Au_\pm)^0 \pm \sum _{p=0}^{k-1} \psi_p^\pm \delta^{(p)}(x_n).
&(1)\cr}$$
The section $\psi_\pm=(\psi_\pm^0,\ldots,\psi_\pm^{k-1})$ of ${\cal L}$
is related to the section $\phi_\pm$ by
$$\psi_\pm={\cal A}\phi_\pm$$
where ${\cal A}$ is a differential operator acting on sections of ${\cal L}$,
and it has block representation
$${\cal A}=\pmatrix{{\cal A}_1&{\cal A}_2&\ldots&{\cal A}_{k-1}&{\cal A}_k\cr
                    {\cal A}_2&{\cal A}_3&\ldots&{\cal A}_k&0\cr
                     \vdots&\vdots&\ddots&\vdots&\vdots\cr
                    {\cal A}_k&0&\ldots&0&0\cr}.\eqno(2)$$
Clearly, the operator ${\cal A}$ is invertible; its inverse is of the form
$${\cal A}^{-1}=\pmatrix{0&\ldots&0&{\cal A}_k^{-1}\cr
 0&\ldots&{\cal A}_k^{-1}& -{\cal A}_k^{-1}{\cal A}_{k-1}{\cal A}_k^{-1}\cr
 \vdots&\ddots&\vdots&\vdots\cr}.$$
We recall that ${\cal A}_k$ is a smooth family of automorphisms.

Now we are ready to prove ${\cal L}_+ +{\cal L}_-={\cal L}$.
Take a section $\phi=(\phi_0,\ldots,\phi_{k-1})$ from ${\cal L}$.
Let $\psi=(\psi_0,\ldots,\psi_{k-1})={\cal A}\phi$, and let $\Psi$
be a distribution
$$\Psi=\sum_{p=0}^{k-1}\psi_p\delta^{(p)}(x_n).\eqno(3)$$
Set $u=A^{-1}\Psi$, and let $u_\pm$ be the restriction of $u$ to $M_\pm$.
It follows from the elliptic regularity theory that $u_\pm$ are smooth
sections of $E_\pm$, up to $\Gamma$. Clearly, $Au_\pm=0$ in $M_\pm$.
Denote by $\phi_j^\pm$ the restriction
of $\partial_n^j u_\pm$ to $\Gamma$. Then sections $\phi^\pm
=(\phi_1^\pm,\ldots,\phi_{k-1}^\pm)$  belong to ${\cal L}_\pm$,
and (1) implies that
$$\psi={\cal A}\phi^+ -{\cal A}\phi^-.$$
Because of invertibility of ${\cal A}$, we conclude that
$\phi=\phi^+-\phi^-.$\hfill\break
$\bigcirc$

Instead of the space ${\cal L}$ of smooth sections of ${\cal E}$,
one can take a Hilbert space
$${\cal H}=H^{k-1+\alpha}(\Gamma,E_\Gamma)\oplus\cdots\oplus
H^{1+\alpha}(\Gamma, E_\Gamma)\oplus H^\alpha (\Gamma,E_\Gamma)$$
of sections. Here $\alpha$ is a sufficiently large positive number.
Let ${\cal H}_\pm$ be the closure of ${\cal L}_\pm$ in ${\cal E}$.
It follows immediately from Proposition 1 and from the standard
elliptic estimates that
$${\cal H}_+\cap{\cal H}_-=\{0\}\quad\hbox{\rm and}\quad
  {\cal H}_+ +{\cal H}_-={\cal H}.$$
To define Sobolev spaces of sections of the vector bundle $E_\Gamma$,
one needs some additional structures: a Riemannian metric on
$\Gamma$, a Hermitian structure on $E_\Gamma$, and a connection
$\nabla^\Gamma$ on $E_\Gamma$. A Riemannian metric and a Hermitian structure
give rise to the $L^2$-scalar product on both $C^\infty(\Gamma, E_\Gamma)$
and $\Lambda^1(\Gamma, E_\Gamma)$, the space of sections of $E_\Gamma$
and the space of one-forms with values in $E_\Gamma$. In the usual way, a
connection induces the operator
$$d^\nabla: C^\infty(\Gamma, E_\Gamma)\to \Lambda^1(\Gamma, E_\Gamma)$$
by the formula
$$d^\nabla=\sum \nabla^\Gamma_{x_i} dx_i$$
in local coordinates. The Laplacian on the space of sections of
$E_\Gamma$ can be defined as $\Delta=(d^\nabla)^* d^\nabla$.
Then, the $H^s$-scalar product is defined as
$$(\phi,\psi)_s=((\Delta+1)^{s}u,v)_{L^2}.$$
Of course, the actual formula for the scalar product depends on all choices
made but it is a well-known fact that spaces themselves are independent
of these choices.

{\bf Proposition 2.} {\sl Let $R_\pm$ be the projection onto} ${\cal H}_\pm$
{\sl parallel to} ${\cal H}_\mp$. {\sl The operators $R_\pm$ are 
pseudodifferential operators, and their complete symbols depend only on
coefficients of $A$ and on their derivatives on the hypersurface $\Gamma$.}

{\bf Proof.} We will treat the projector $R_+$; clearly, the case of
$R_-$ is similar. Let us define a new operator $\tilde R$ acting
on sections of ${\cal E}$. We describe now how to construct
$\zeta=(\zeta_0,\ldots,\zeta_{k-1})=\tilde R\phi$ where
$\phi=(\phi_0,\ldots,\phi_{k-1})$. Firstly, we define
$\psi={\cal A}\phi$ (see (2) for ${\cal A}$), then the section $\psi$
is used to produce the distribution $\Psi$ (see (3)), then the section
$u$ of $E$ is defined by $u=A^{-1}\Psi$, and, finally, $\zeta_j$
is the restriction to $\Gamma$ of $\partial_n^j u_+$ where
$u_\pm$ are the restrictions of $u$ to $M_\pm$. We will show that,
in fact, $\tilde R=R_+$, and then we will see that $\tilde R$ is a
pseudodifferential operator, and we will discuss how to compute
its symbol.

To show that $\tilde R=R_+$ one has to verify that $\tilde R\phi=\phi$
when $\phi\in{\cal H}_+$ and $\tilde R\phi =0$ when $\phi\in {\cal H}_-$.
Let $\phi\in {\cal H}_+$. Then there exists a solution $v$ of the equation
$Av=0$ in $M_+$ such that $\phi_j$ is the restriction to $\Gamma$
of $\partial_n^j v$. Let $v^0$ be the section of $E$ that equals $v$
over $M_+$, and that equals $0$ over $M_-$. Then $Av^0=\Psi$.
Hence, $u=v^0$, $u_+=v$, and $\zeta_j=\phi_j$. 

Now, let $\phi\in{\cal H}_-$. Then there exists a solution $w$ of the
equation $Aw=0$ in $M_-$ such that $\phi_j$ equals the restriction
of $-\partial_n^j w$ to $\Gamma$. Let $w^0$ be the section of $E$
that coincides with $w$ over $M_-$, and that vanishes over $M_+$.
Then $Aw^0=\Psi$, so $w=u$, and $u_+=0$. We conclude that 
$\zeta=\tilde R\phi=0$.

Denote by $r^+$ the operator of restricting a section from $M_+$ to
$\Gamma$. Let $B_{qp}$ be the operator acting on sections of $E_\Gamma$
according to the formula
$$B_{qp}\psi_p=r^+\partial_n^q A^{-1}(\psi_p\delta^{(p)}(x_n)).$$
Let ${\cal B}$ be the operator acting on sections of ${\cal E}$,
the block-matrix of which is $(B_{qp})$. Clearly,
$$R_+=\tilde R ={\cal B}{\cal A}.\eqno(4)$$
Denote by $S(x',x_n;\xi',\xi_n)$ the symbol of $A^{-1}$ in certain
local coordinate system, with respect to a local frame.
Up to a smoothing operator applied to $\psi_p$,
$r+\partial_n^q A^{-1}(\psi_p(x')\delta^{(p)}(x_n)$ equals
$$\lim_{x_n\to 0^+}(2\pi)^{-n}i^{p+q}\int\xi_n^{p+q}e^{i(x'-y')\xi'}
e^{xi_n\xi_n}S(x',x_n;\xi',\xi_n)dy'd\xi'd\xi_n.$$
It is a well known fact that such an operator is a pseudodifferential
operator (e.g. see [10]), and its symbol equals
$$\lim_{x_n\to 0^+}{i^{p+q}\over 2\pi}\int_{-\infty}^\infty
\xi_n^{p+q}e^{ix_n\xi_n}S(x',x_n;\xi',\xi_n)d\xi_n.\eqno(5)$$
The complete asymptotic symbol $S$ is meromorphic in $\xi_n$,
and the expression (5) can be rewritten as
$$\sigma(B_{qp})={i^{p+q}\over 2\pi}\int_{\gamma_+}S(x',0;\xi',\xi_n)
                 d\xi_n \eqno(6)$$
where $\sigma$ means ``symbol'', and $\gamma_+$ is a contour in
the complex $\xi_n$-plane that goes around all poles of $S$ in the
counter-clockwise direction. The symbol of ${\cal B}$ depends only
on the restriction of the symbol of $A^{-1}$ to $\Gamma$, and this restriction
depends only on coefficients of $A$ and on their derivatives on $\Gamma$.
\hfill\break
$\bigcirc$

Matrix elements, $(R_+)_{qj}$, of the operator $R_+$ are given by
$$(R_+)_{qj}=\sum_{p=0}^{k-1-j}B_{qp}{\cal A}_{p+j+1}\eqno(7)$$
(cf. (2) and (4)).
It can be easily seen from (6) or (5) that the order of $B_{qp}$
equals $p+q-k+1$. The order of the operator ${\cal A}_{p+j+1}$
equals $k-p-j-1$, so each term on the right in (7) is a pseudodifferential
operator of order $q-j$. We conclude that
$$\hbox{\rm ord} (R_+)_{qj}=q-j.$$
It is clear from (2) and (6) that the principal symbol of
$(R_+)_{qj}$ depends only on the value of the principal symbol of $A$ on 
$\Gamma$.

If a pseudodifferential operator has a block-matrix form,
and the order of its $(qj)$-entry is $q-j$ then by a principal
symbol of such an operator we mean a matrix, the $(qj)$th entry
of which is the $(q-j)$-principal symbol of the corresponding
operator (this is the principal symbol in the sense of
Agmon--Douglis--Nirenberg). Such an operator is similar to an operator
of order $0$.

Now, let $P_\pm$ be the orthogonal projector onto the space ${\cal H}_\pm$.

{\bf Theorem 1.} $P_\pm$ {\sl is a pseudodifferential operator; its $(qj)$-th
entry in the block-matrix representation has order $q-j$. The complete
symbol of $P_\pm$ depends only on the coefficients of $A$ and 
on their derivatives
on $\Gamma$, and the principal symbol of $P_\pm$ depends only
on the value of the principal symbol of $A$ on $\Gamma$.}

{\bf Proof.} Because operators $R_\pm$ are projectors, their principal symbols
are also projectors. Therefore, for any value of $\lambda\not=0,1$,
the operator $\lambda-R_\pm$ is elliptic in the sense of 
Agmon--Douglas--Nirenberg, and the resolvent $(\lambda-R_\pm)^{-1}$
is a holomorphic in $\lambda$ family of pseudodifferential operators 
that have the same structure. All statements of the theorem follow
from the formula
$$P_\pm={1\over 2\pi i}\int_{\gamma}(\lambda-R_\pm)^{-1}d\lambda$$
where $\gamma$ is, say, a circle of radius $1/2$ centered at the origin
and taken with the counter-clockwise orientation.
\hfill\break
$\bigcirc$
\vfill\eject

\centerline{\bf 3. A manifold with boundary}
Now, instead of treating a closed manifold that is divided into two
connected components by a hypersurface, we will consider a
compact manifold $M_0$ with smooth boundary $\Gamma$.
Let $E_0$ be a vector bundle over $M_0$, let $A$ be an elliptic differential
operator of order $k$ acting on sections of $E$
that satisfies the Agmon--Seeley condition. By ${\cal L}_0$ we denote
the subspace of ${\cal L}$ that consists of sections
$(u,\partial_n u,\ldots,\partial_n^{k-1}u)$ where $Au=0$ in $M_0$,
and $P_0$ is the orthogonal projection onto ${\cal L}_0$.
All notations are the same as those used above, with the only
difference that a neighborhood of $\Gamma$ is diffeomorphic to
$\Gamma\times [0,1)$ (not to $\Gamma\times (-1,1)$), and we use
the subscript $0$ instead of $+$ or $-$. We will prove the following
theorem.

{\bf Theorem 2.} $P_0$ {\sl is a pseudodifferential operator; its $(qj)$-th
entry in the block-matrix representation has order $q-j$. The complete
symbol of $P_0$ depends only on the coefficients of $A$ and 
on their derivatives
on $\Gamma$, and the principal symbol of $P_0$ depends only
on the value of the principal symbol of $A$ on $\Gamma$.}

{\bf Proof.} We will construct a closed manifold $M\supset M_0$,
a vector bundle $E\to M$ that extends $E_0$, and an extension of $A$
to an invertible, selfadjoint, elliptic operator acting on sections
of $E$. Then the statements of Theorem 2 will follow from Theorem 1
(just replace the subscript $0$ by $+$).

For the manifold $M$, we take the double of $M_0$, and, for the vector
bundle $E$, we take the double of $E_0$. A neighborhood
of $\Gamma$ in $M$ is diffeomorphic $\Gamma\times (-1,1)$.
Clearly, $A$ can be extended
to en elliptic differential operator that satisfies the Agmon--Seeley
condition, and that acts on sections of $E$ over
$M_0\cup\Gamma\times (-1,0)$. We will construct an extension
of $A$ to the whole manifold $M$ in two steps.

{\bf Step 1.} Let $a(x,\xi)$ be the principal symbol of $A$,
and let $\gamma$ be a positively oriented contour in the complex plane
that does not intersect a ray $z=re^{i\theta}$, $r\geq 0$, and that
encloses all eigenvalues of $a(x,\xi)$ for all $x$ and all $\xi$,
$|\xi|=1$. Here, we have chosen an arbitrary metric on the cotangent
bundle to $M$. The existence of such a contour
is guaranteed by the Agmon--Seeley condition.
Let $\chi(\tau)$ be a smooth non-negative function of one variable that equals
$0$ when $\tau<-3/4$, and that equals $1$ when $\tau>-1/2$.
We define a symbol $b(x,\xi)$ on the cosphere bundle $S^*M$
(it is an endomorphism of the vector bundle $E$, pulled back to $S^*M$)
in the following way: $b(x,\xi)=a(x,\xi)$ over $M_0$, $b(x,\xi)$ equals
the identity operator over $M\setminus (M_0\cup\Gamma\times(-1,0))$, and
$$b(x,\xi)={1\over 2\pi i}\int_\gamma z^{\chi(x_n)} (z-a(x,\xi))^{-1}dz$$
over $\Gamma\times (-1,0)$. To define the powers of $z$, one makes
the cut $re^{i\theta}: r\geq 0$ in the complex plane.
The symbol $b(x,\xi)$ is extended to the whole cotangent bundle
by $k$-homogeneity. Let $B$ be a pseudo-differential operator
with the principal symbol $b(x,\xi)$. Clearly, the operator
$B$ is elliptic, it satisfies the Agmon--Seeley condition, and its
principal symbol equals $a(x,\xi)$ over $M_0\cup\Gamma\times (-1/2,0)$.

{\bf Step 2.} Choose smooth function $\phi(x)$ and $\psi(x)$ on $M$
such that $\phi(x)=1$ on $M_0\cup\Gamma\times (-1/4,0)$,
$\phi(x)=0$ outside of $M_0\cup\Gamma\times (-1/2,0)$, and
$\phi^2(x)+\psi^2(x)=1$. Define the operator 
$$\tilde A=\phi(x)A\phi(x)+\psi(x)B\psi(x).$$
The operator $\tilde A$ is a pseudo-differential elliptic operator
that satisfies the Agmon--Seeley condition, it is differential
in $\Gamma\times (-1/4,1)$, and it coincides with $A$ over $M_0$.
If it is invertible, then one can apply Theorem 1.
If it is not invertible, then one can make it invertible by adding
to it an operator of multiplying by an appropriate function
supported on $M\setminus M_0$.\n
$\bigcirc$

{\bf Remark 1.} The Theorem 2 holds for pseudo-differential operators
that are differential in a neighborhood  of $\Gamma$. In fact,
in the proof we did not use the fact that the operator $A$ 
is differential outside $\Gamma\times [0,1)$.

{\bf Remark 2.} The statements of the Theorem 1 depend only
on the restrictions of $A$ to $M_\pm$. Theorem 2 shows (see the previous
remark) that, for conclusions of Theorem 1, invertibility of $A$
is not essential. If the operator $A$ is not invertible, then,
in Proposition 1, ${\cal L}_+\cap{\cal L}_-$ is finite-dimensional,
and ${\cal L}_+ +{\cal L}_-$ has finite codimension.

{\bf Remark 3.} One can remove the assumption that the operator $A$
is self-adjoint. In fact, the operator
$$A'=\pmatrix{0&A\cr
              A^*&0\cr}$$
is a self-adjoint operator acting on sections of the vector bundle $E\oplus E$.
The statements of the Theorem 2, when applied to the operator $A'$,
imply the same statements for the operator $A$.

Now, let us go back to the situation when there are two differential
operators around, $A$ and $\tilde A$; they act on sections
of vector bundles $p_E:E\to M$ and $p_{\tilde E}\tilde E\to \tilde M$, and
$\partial M=\partial\tilde M=\Gamma$. We assume that there exists an
isomorphism $\Psi: E_U\to \tilde E_{\tilde U}$ such that
$p_E\circ\Psi=p_{\tilde E}$ over $\Gamma$. Here $U$ is a neighborhood
of $\Gamma$ in $M$ and $\tilde U$ is a neighborhood of $\Gamma$
in $\tilde M$. We also assume that the isomorphism $\Psi$ maps
the connection $\tilde\nabla$ on $\tilde E$ to the connection
$\nabla$ on $E$. To avoid confusion, let us say that in the general setting
we use connections for completely different purposes than they were
used in the context of the Dirac operator. For the Dirac operator,
the connection was used to construct it. Here operators are given
by their symbols, and the connection is used for the purpose
of taking normal derivatives on the boundary. The isomorphism
$\Psi$ identifies $U$ and $\tilde U$, $E_U$ and $\tilde E_{\tilde U}$,
so we will think of $E_U$ and $\tilde E_{\tilde U}$ as being the same.

Let
$$a(x,\xi)=a_k(x,\xi)+a_{k-1}(x,\xi)+\cdots+a_0(x,\xi)$$
be the splitting of the total (complete) symbol of the operator $A$ into
homogeneous components in certain local coordinates. In this splitting,
only the principal symbol $a_k$ has an invariant meaning.
We say that two operators $A$ and $\tilde A$ {\it agree up to the
order $q$ on the boundary} $\Gamma$, if they have the same order, and
the corresponding homogeneous components,
$a_{k-j}$ and $\tilde a_{k-j}$, are equal on $\Gamma$, together
with all their derivatives up to the order $q-j$, for $j\leq\min\{k,q\}$.
Though the components $a_{k-j}$ themselves are not defined invariantly,
the property of two operators to agree up to a certain order
on $\Gamma$ does not depend on a particular choice of local coordinates.
Note that Dirac operators that were discussed in the introduction
agree up to the order $0$ on $\Gamma$.

Denote by $H_+^A$ the closure in ${\cal H}$ of the space of Cauchy data
for the operator $A$ (of the space ${\cal L}_0$ from the Theorem 2),
and let $H_+^{\tilde A}$ be the closure of the space of Cauchy data
for the operator $\tilde A$.
To measure, how far the space $H_+^A$ is from the space
$H_+^{\tilde A}$, we introduce orthogonal projections
$\pi_+^A$ and $\pi_+^{\tilde A}$ onto these spaces. It follows from
Theorem 2 that\n
{\sl if the operators $A$ and $\tilde A$ agree up to the order
$q\geq 0$ on $\Gamma$ then $\pi_+^A-\pi_+^{\tilde A}$ is a compact operator,
and it belongs to the Schatten ideal $\Sigma_p$ for $p>(n-1)/(q+1)$.}

>From the construction of projections $R_\pm$ and $P_\pm$ in section 2
it follows that the first $q+1$ terms in their complete symbols depend
on the restriction of $a_{k-j}(x,\xi)$, and its derivatives up
to the order $q-j$, to $\Gamma$ where $j\leq\min\{k,q\}$.
The same is true for the operator $P_0$ from this section.
It follows that if operators $A$ and $\tilde A$ agree up to
the order $q$ on $\Gamma$ then $\pi_+^A-\pi_+^{\tilde A}$ is a 
pseudo-differential
operator, and the $(ij)$-th entry in its block matrix representation
has the order $i-j-q-1$. It is well known that the singular numbers $s_j$ 
of such
an operator can be estimated
$$s_j\leq Cj^{-(q+1)/(n-1)}$$
(note that the dimension of $\Gamma$ equals $n-1$),
and, therefore, it belongs to $\Sigma_p$ for $p>(n-1)/(q+1)$.

Let us now denote by $P_{\tilde A A}$ the restriction of the projection
$\pi_+^A$ to $H_+^{\tilde A}$ and by $Q_{\tilde A A}$
the restriction of $I-\pi_+^A$ to $H_+^{\tilde A}$. Then, the operator
$P_{\tilde A A}$ is Fredholm. In fact, its kernel coincides with
$H_+^{\tilde A}\cap (H_+^A)^\perp$. The compact operator 
$\pi_+^{\tilde A}-\pi_+^A$
is identical on this space; therefore it is finite-dimensional.
The adjoint to $P_{\tilde A A}$ is $P_{A \tilde A}$, and, by the same
reason, its kernel is finite dimensional. 
Further, if 
$\pi_+^A-\pi_+^{\tilde A}$
belongs to the Schatten class $\Sigma_p$ then $Q_{\tilde A A}\in \Sigma_p$.
In fact, $Q_{\tilde A A}$ is the restriction to $H_+^{\tilde A}$ 
of the operator
$(I-\pi_+^A)\pi_+^{\tilde A}=(\pi_+^{\tilde A}-\pi_+^A)\pi_+^{\tilde A}$ 
which belongs to
$\Sigma_p$ because $\Sigma_p$ is an ideal in the ring of bounded operators.

{\bf Remark 4.} The index of the operator $P_{\tilde A A}$ needs not
be equal to $0$. On the other hand, the index of the corresponding projections
constructed for the operator $A'$ (see Remark 3) always equals $0$.

\centerline{\bf References}
\n
[1] Yu. Manin, Funk. Anal., {\bf 20} (1986), p. 88\n
[2] C.Vafa, Phys. Lett., {\bf B190}, p. 47\n
[3] L.Alvarez--Gaume, C.Gomez, C.Reina, Phys. Lett., {\bf B190} (1987), p. 55\n
[4] A.Morozov, JETP Lett., {\bf 45} (1987), p. 585\n
[5] A.Schwarz, JETP Lett., {\bf 46} (1987), p. 438\n
[6] A.Schwarz, Grassmannian and String Theory, preprint, hep-th., 9610122\n
[7] J.Mickelsson, Comm. Math. Phys, {\bf 127} (1990), p. 285\n
[8] M.Atiyah, V.Patodi, I.Singer, Math. Proc. Camb. Phyl, Soc.,
   {\bf 77} (1975), p. 43\n
[9] R.T.Seeley, Amer. J. Math., {\bf 88} (1966), p. 781\n
[10] S.Rempel, B.-W.Schulze, {\it Index theory of elliptic boundary problems},
Academie--Verlag, Berlin, 1982\n

\vfill\eject\end